\documentclass[11pt, letterpaper]{amsart}
\usepackage{amsmath}
\usepackage{amssymb}
\usepackage{color}

\newcommand{\qqed}[0]{ \rule{6pt}{8pt}}
\newcommand{\cdott}[0]{\cdot\cdot\cdot}
\newcommand{\lrar}[0]{\longrightarrow}

\newcommand{\sub}[0]{\subseteq}
\newcommand{\sm}[0]{\backslash}
\newcommand{\lb}[0]{\linebreak}
\newcommand{\prr}{\partial} 
\newcommand{\krr}{\textrm{Ker}} 
\newcommand{\iim}{\textrm{Im}}  

\newcommand{\ho}{\hspace{1pt}}
\newcommand{\dyi}{\hspace{-2pt}:=}
\newcommand{\lss}{\vspace{1mm}}
\newcommand{\cbk}{\color{black}}
\newcommand{\cbl}{\color{blue}}

\begin{document}
\setlength{\lineskiplimit}{5pt}
\setlength{\lineskip}{5pt}
\centerline{\large\bf The Permanent Rank of a Matrix (Part Three)}

\centerline{\large\bf Note on the Additive Basis Conjecture} 

\bigskip
\centerline{Yang\, Yu\footnote{contact: yang.yu30@rutgers.edu, Dept of Math, Rutgers University}
\hspace{11mm}Date\; 2026-05-27}

\bigskip
{\bf Abstract.}\; We show that in a vector space over $Z_3$, the union of any four linear bases 
is an additive basis, thus proving the Additive Basis Conjecture for $Z_3$, 
and providing an alternative proof of the weak 3-flow conjecture.

\bigskip\medskip
\centerline{\large\bf  Introduction}

\bigskip\medskip
The Additive Basis Problem is a classical problem in additive combinatorics whose history parallels
that of the more famous Arithmetic Progressions Problem. Both have been extensively studied since the
1930's (see e.g. [Erdös]), first for the integers, but later also for other abelian groups; 
see e.g. [Mesh] (for arithmetic progressions) and [JLPT] (for additive bases) for early work taking 
this more general viewpoint.

For arithmetic progressions, recent years have seen celebrated progress on the central problem of 
bounding the sizes of AP-free sets in $Z_p^n$  [CLP, EG]. But there has been no comparable 
breakthrough on what's perhaps the best-known problem on additive bases in $Z_p^n$: the following 
conjecture of Jaeger, Linial, Payan and Tarsi.

Recall that a multiset $B$ is an {\em additive basis}  of a vector space $S$, if every element of $S$
is a linear combination of elements in $B$, with each coefficient either 0 or 1.

\cbl{\bf Conjecture 1}\quad (The Additive Basis Conjecture [JLPT])

For any prime $p$, there exists a constant $c(p)$, such that in any vector space
over $Z_p$,  the multiset union of any $c(p)$ linear bases is an additive basis.\cbk

Conjecture 1 was studied in [ALM, Sz, NPT1, NPT2, EVLT, HQ, CKMS, Ma, LMP]. It is related to a few other 
problems in discrete mathematics. For example, the case $p=3$ implies Francois Jaeger's
weak 3 Flow Conjecture (first proved by Carsten Thomassen in 2012 [Th] by graph theory). 

It is proved in  \cite{a} that the union of any $c(p)\log n$ linear bases is an additive basis,
where $n$ is the dimension of the vector space.
Our approach, like that of \cite{a}, is based on permanents. 
Define the {\em perrank} of a matrix $M_{m\times n}$ to be the size of a largest 
square submatrix with 
nonzero permanent; if this is equal to $m$ or $n$, then we say $M$ has {\em full perrank}.

\cbl{\bf Conjecture 2}\; \cite{a}\; If\, $B_1,B_2,\cdott, B_p$ are nonsingular matrices over a field of 
characteristic $p\ge 3$,\, then $\left(\begin{array}{cccc} 
B_1 & B_2 &\cdott & B_p \\ \vdots &\vdots &&\vdots \\ B_1 & B_2 &\cdott & B_p\end{array}\right)$,
 where each row repeats $p-1$ times, has full perrank.\cbk

\lss
Here we give the first constant bounds for both conjectures, and in particular the first proof of {\em any}
case of the Additive Basis Conjecture:

\lss
\cbl{\bf Theorem 3\; (Main Theorem)}\quad If\, $P,R,S,T$ are nonsingular matrices over a field of 
characteristic 3, then $\left(\begin{array}{cccc} P & R & S & T \\ P & R & S & T \end{array}\right)$
has full perrank.\cbk

\lss
Theorem 3 gives Conjecture 1 when $p=3$ with $c(3)=4$, by the Combinatorial Nullstellensatz 
(see \cite{a} section 3 for details), while Conjecture 2 implies Conjecture 1 with $c(p)=p$. 
Theorem 3 will be proved at the end of the paper, after we have developed the necessary machinery,
the main point here being Theorem 7. An early look at the easy derivation of Theorem 3 should help
to motivate what precedes it.

This paper is part 3 of the author's series ``The permanent rank of a matrix”;
part 1 was [Yu]; and at this writing part 2 is still in preparation.
 
\bigskip\medskip
\centerline{\large\bf  Definitions and Notation}
\bigskip\medskip

Given a field, let $A^n$ be the quotient of the polynomial ring in $n$ variables \lb
$x_1, x_2, \cdott, x_n$  by the ideal generated by $x_1^2,\,x_2^2, \cdott, x_n^2$.  
The $k$\,-\,th degree component of the graded algebra $A^n$  is denoted by $A_k^n$;\, 
we omit $n$ when there is no ambiguity.
For an ideal $J\sub A$, let $J_k=A_k \cap J$.

For $f\in A$, define $\krr(f)=\{g:gf=0\}$,\; $\iim(f)=\{ fg: g\in A\}$. 

\lss
We introduce two operators. Let  $\prr_x$ and $E_x$ be the quotient 
and remainder of formal  division by $x$; we use $\prr_i$ for $\prr_{x_i}$ and similarly for $E_i$. 
For example, if $f=x_1x_2+x_1x_3+x_2x_3$,  then $\prr_1 f=x_2+x_3$\, and\, $E_1 f=x_2x_3$.  

(We use the letter $E$ because $E_i$ {\em eliminates}  all terms containing $x_i$).

\lss
It is obvious that\; $E_i(fg)=(E_i f)(E_i g),\;\;\prr_i(fg)=(\prr_i f)(E_i g)+(\prr_i g)(E_i f)$, and
$E_i E_j=E_j E_i,\; \prr_i\prr_j=\prr_j\prr_i,\;   E_i\prr_j=\prr_j E_i$\;
(but $E_i \prr_j \neq E_j \prr_i$).

For $u\in A_1$, define its {\em support}  to be supp$(u)\dyi\{ x: \prr_x u\neq 0\}$.
An element of $A_1$ is also called a {\em linear form}.

Let $u$ be a linear form with  $\prr_x u=c\neq 0$, set $\prr_x=\prr,  E_x=E$ for simplicity. 
For any $f$, define another division operation: divide $f$ by  $u$  w.r.t. $x$ by

$f=Ef+(\prr f) x=Ef+(\prr f)(u-Eu) c^{-1}=c^{-1}(\prr f) u + Ef-c^{-1}(\prr f)(Eu)$\;
and define  $R_{(u,\,x)} f\dyi Ef-c^{-1}(\prr f)(Eu)$  as the remainder. 

Evidently\, $\prr(R f)=0$  and $f-Rf\in\iim(u)$, this is what we need later.

\lss
For $U$ and $V$ subspaces of $A_i$ and $A_j$ respectively, define $UV$ to be the
subspace of $A_{i+j}$ spanned by $\{uv: u\in U,\;v\in V\}$. Define 
$\iim(U)$ to be the ideal generated by $U$, and $\krr(U)=\{f:fu=0\;\;\forall\,u\in U\}$.

\lss
A subspace of $A_1$ is also called a {\em linear form space}. For a linear form space $U$, 
define its {\em support}  to be supp$(U)\dyi\hspace{-1mm}\bigcup\limits_{u\in U}$supp$(u)$;
\vspace{-2mm} define its {\em minimum support function} to be 

ms$\ho_i(U)\dyi$ min\,$\{|\,\textrm{supp}(V)|: V\sub U\;\textrm{with dim}(V)=i\,\}$
for $i\le $ dim$(U)$.

\lss
Label the rows and columns of a matrix $M_{m\times n}$ with variables $x_1,x_2,\cdott,x_m$ and
$y_1,y_2,\cdott,y_n$ respectively, and view its rows and columns as linear forms in $A^n$ and $A^m$. 
Then $M$ has full perrank  iff the product of its rows or columns is nonzero in  $A^n$ or $A^m$.

\medskip
Since direct induction does not work for Conjecture 2, we attempt to induct on a stronger 
or equivalent statement that is suitable for induction. Specifically, we use the minimum 
support function and matroid theory to formulate an equivalent statement of Conjecture 2.

\cbl{\bf Conjecture 4}\;\; Suppose $U$ is a linear form space over a field of characteristic 
$p\ge 3$, dim$(U)=n$, and ms$_i(U)\ge p i$  for all $1\le i\le n$, then $U^{(p-1)n}\neq 0$.\cbk

One direction, Conjecture 4 implies Conjecture 2:
Let $U$ be the linear form space spanned  by the rows of $(B_1,\, B_2,\cdott, B_p)$, then $U$
satisfies the condition of Conjecture 4, and $U^{(p-1)n}\neq 0$ gives full perrank.

The other direction: Let $S=\text{supp}(U)$, choose a linear basis of  $U$, 
we get a matrix $M_{n\times S}$. There is a matroid structure on $S$ with rank function $r()$.

For any $R\sub S$, let $r=r(R)$, which is the rank of the $n\times R$ submatrix. Perform 
row operations on $M$ to make the bottom  $n-r$ rows of the $n\times R$ submatrix zero.
Let $V\sub U$  be the subspace spanned by the bottom $n-r$ rows of $M$. We have
$|S\ho\sm R| \ge |\text{supp}(V)| \ge p(n-r) = p(n-r(R))$.

By Edmond matroid base packing theorem, $S$ contains $p$ disjoint bases.
Let $B_1, B_2,\cdott, B_p$ be the nonsingular matrices corresponding to these bases,
then Conjecture 2 gives $U^{(p-1)n}\neq 0$.\qqed

\lss
The equivalence is the key idea of this paper, as well as the motivation of Theorem 7.

\newpage
\centerline{\large\bf  Supporting Results and Proofs}
\bigskip\medskip

From now on, we assume the ground field has characteristic 3 and is infinite,
otherwise extend it to infinity.  The following theorem is the base step for the main induction 
in the proof of Theorem 7.

\lss\cbl{\bf Theorem 5}

{\bf (1)}\, $\krr_k(u)=\iim\ho_k(u^2)$\, for any linear form $u$ with $|\,$supp$(u)|\ge 2k+1$.

{\bf (2)}\, $\krr_k(u^2)=\iim\ho_k(u)$\, for any linear form $u$ with $|\,$supp$(u)|\ge 2k+2$.\cbk

\lss{\em Proof.}\; Since $u^3=0$, one direction is trivial. The other direction is
by induction on $k$,\, easy to verify when $k\le 1$.  
Pick any  $x\in$ supp$(u)$, and set $\prr_x=\prr,  E_x=E$ for simplicity. WMA\, $\prr  u=1$.

(1) Suppose $f\in \krr_k(u),\; fu=0$; take $\prr,\; E f+(\prr f)(E u)=0$; multiply by 
$E u,\; (\prr f)(E u)^2=0$.
By induction hypothesis of (2),\, $\prr f=g(E u)$ for some $g$, then\newline
$f=Ef+(\prr f) x=(\prr f)(x-Eu)=-g(Eu)(Eu-x)=-g(Eu+x)^2=-g u^2$.

(2) Suppose $f\in \krr_k(u^2),\; fu^2=0$; take $\prr,\; (2Ef+(\prr f)Eu)Eu=0$. By (1) we have
$Ef-(\prr f)Eu=g(Eu)^2$ for some $g$, then\newline
$f=Ef+(\prr f) x=(\prr f)(Eu+x)+g(Eu)^2=(\prr f)u+g u(Eu-x)\in\iim(u)$.\qqed
 
We say a linear form space $U$ {\em covers}  $(a_1,a_2,\cdott,a_k)$  if\,
ms$\ho_i(U)\ge a_i$ for all $1\le i\le k$. 

The following lemma plays a crucial role in the proof of Theorem 7.

\lss
\cbl {\bf Lemma 6}\quad Suppose $U$ is a linear form space with dim$(U)=n$  covering an
increasing sequence $(a_1,a_2,\cdott,a_n)$.  Then for each $0\le k\le n$, there exists a subspace
$U_k\sub U$ with dim$(U_k)=k$  covering $(a_{n+1-k},  a_{n+2-k}, \cdott, a_n)$.\cbk

\lss
{\em Proof.}\; Set  $U_0=0$ and suppose $U_k$ exists.

For each $S\sub$ supp$(U)$ with $|S|=a_{n-k}-1$, let

$V_S\dyi\{v : v\in U,\textrm{ there exists }\,u\in U_k\,\textrm{ such that supp}(v+u)\sub S\ho\}$;
\newline note $U_k\sub V_S$ since supp$(0)=\varnothing$. 
Claim  $V_S$  is a proper subspace of $U$, otherwise 
choose $\{v_i\}$ such that $U=$ span$(U_k,  v_1, v_2, \cdott, v_{n-k})$.
For each $i$, choose  $u_i\in U_k$  such that supp$(v_i+u_i)\sub S$.
Let $V=$ span$(\{v_i+u_i\})$, then dim$(V)=n-k,\, |\,$supp$(V)|<a_{n-k}$,\, contradiction. 

Because the ground field is infinite, $U$ is not a finite union of its proper subspaces.
Choose   $v\in U$  that is not in any $V_S$,  and 
let  $U_{k+1}\dyi$ span$(U_k,  v)$.

If\, $u\in U_{k+1}\sm U_k, \textrm{then }|$\,supp$(u)|\ge a_{n-k}$  by our choice of  $v$. 
If\, $0\neq u\in U_k$,  then by induction hypothesis,
$|$\,supp$(u)|\ge$ ms$_1(U_k)\ge a_{n+1-k}>a_{n-k}$.  So ms$_1(U_{k+1})\ge a_{n-k}$.

For any $H\sub U_{k+1}$ with dim$(H)=h\ge 2$,  either  dim$(H\cap U_k)=h-1$
or $H\sub U_k$.  By induction hypothesis,  either 

$|$\,supp$(H)|\ge |\,$supp$(H\cap U_k)|\ge$ ms$\ho_{h-1}(U_k)\ge a_{n+h-1-k}$, or

$|$\,supp$(H)|\ge$ ms$\ho_h(U_k)\ge a_{n+h-k}>a_{n+h-k-1}$.

So  ms$\ho_h(U_{k+1})\ge a_{n+h-k-1}$.\qqed

\vspace{3mm}
\cbl{\bf Theorem 7}\quad Suppose $U$ is a linear form space, dim$(U)=n$ and $k\ge0$.

{\bf (A)}\, If\, ms$\ho_i(U)\ge 4 i-2+2k$\, for $1\le i\le n$, then 
$\krr_k(U^{2n})=\iim\ho_k(U)$.\newline\indent
{\bf (B)}\, If\, ms$\ho_i(U)\ge 4 i-3+2k$\, for $1\le i\le n$, then \newline
$\krr_{2n-2+k}(U)=\iim\ho_{2n-2+k}(U^{2n})$.\cbk

\lss
The proof is delicate, any mismatch between degree and support or other discrepancy
invalidates it. We need to check degree and support\,(abbrev. CDS) 
7 times; 5 times exact match; twice there is extra support of exactly one. 
The induction hypothesis is applied 6 times in the proof.

\lss
{\em Proof.}\; One direction is trivial, the other direction is by induction on $n$.
Theorem 5 gives the case $n=1$. Suppose true for $n-1$, then induction on 
$k$ by : $B(n,0)\lrar A(n,0)\lrar B(n,1)\lrar A(n,1)\lrar B(n,2)\lrar A(n,2)
\cdott\cdot$

\lss\cbl
Claim: $A(n,k-1)$ implies $B(n,k)$ for all $k\ge 1$.\cbk

Suppose $U$ satisfies condition (B) and $f\in\krr_{2n-2+k}(U)$. By Lemma 6, there
exists $V\sub U$ with dim$(V)=n-1$, and ms$\ho_i(V)\ge 4 i+1+2k$\, for all $1\le i\le n-1$.
Choose a variable $x$ and a linear basis $\{u_1+x, u_2,\cdott, u_{n-1}\}$ of $V$ 
such that $\prr_x u_i=0$  for all $1\le i\le n-1$. Choose any  $u_n\in U\sm V$ with
$\prr_x u_n=0$. Then $\{u_1+x, u_2,\cdott, u_n\}$ is a linear basis of  $U$.

Note\, $2n-2+k=2(n-1)-2+(k+2)$, and $f\in\krr_{2n-2+k}(V)$ also. 
Observe\, $V^{2(n-1)}=$ span$((u_1+x)^2 u_2^2 \cdott u_{n-1}^2)$.
Apply $B(n-1,k+2)$ to $V$, CDS exact match. We get 
$f=g(u_1+x)^2 u_2^2 \cdott u_{n-1}^2$ \hfill (1)\newline
for some $g$ with deg$(g)=k$. WMA\, $\prr_x g=0$, otherwise replace it with $R_{(u_1+x,\,x)}g$.
Then $fu_n=0$  gives\, $g u_n (u_1+x)^2 u_2^2 \cdott u_{n-1}^2=0$.

Apply $A(n-1,k+1)$ to $V$, CDS  extra support of one. We have\newline
$gu_n=a_1(u_1+x)+a_2 u_2+\cdott+a_{n-1}u_{n-1}$\; for some $a_i$.

\lss\cbl If\, $k=0$, then $g=a_i=0,\; f=0$. Here we got $B(n,0)$.\cbk 

\lss
If\, $k\ge 1$, take $E_x$ and $\prr_x$\,(write $E_x, \prr_x$ as $E, \prr$), we get 

$g u_n=E(a_1)u_1+E(a_2)u_2+\cdott+E(a_{n-1})u_{n-1}$;\; and

$0=E(a_1) +\prr(a_1)u_1+ \prr(a_2)u_2+\cdott+\prr(a_{n-1})u_{n-1}$,\; multiply by $u_1$, 
\newline
we get\, $g u_n+\prr(a_1)u_1^2 \in \iim($span$(u_2, \cdott, u_{n-1}))$\hfill (2)

Multiply by $u_2^2\cdott u_{n-1}^2 u_n^2$, we get\; $\prr(a_1) u_1^2 u_2^2\cdott u_n^2=0$.

Since ms$\ho_1(U)\ge 1+2k\ge 3$, $x\notin U$, so dim$(E_x(U))=n$.
Apply $A(n,k-1)$ to $E_x(U)$, CDS exact match. We get\,
$\prr(a_1)=c_1 u_1+c_2 u_2+\cdott+c_n u_n$.

Substitue into (2), we get $(g+c_n u_1^2)u_n\in\iim($span$(u_2, \cdott, u_{n-1}))$.

Multiply by $u_2^2\cdott u^2_{n-1}$, we get 
$(g+c_n u_1^2)u_2^2\cdott u_{n-1}^2 u_n=0$.  Introduce a dummy 
variable $y$ to make $(g+c_n u_1^2) u_2^2\cdott u_{n-1}^2 (u_n+y)^2=0$.

Let $Y\dyi$ span$(u_2,\cdott, u_{n-1}, u_n+y)$. For any $I\sub Y$ with dim$(I)=i$,
if\, $y\notin$ supp$(I)$, then $I\sub$ span$(u_2,\cdott, u_{n-1})\sub V$ with 
$|\,$supp$(I)|\ge 4 i+1+2k$. If\, $y\in$ supp$(I)$, then 
$|\,$supp$(E_y(I))|\ge 4 i-3+2k$\, since $E_y(I)\sub U$, and so $|\,$supp$(I)|\ge 4 i-2+2k$.
Apply $A(n-1,k)$ to $Y$, CDS exact match. We have\, $g+c_n u_1^2\in\iim(Y)$. Take $E_y$, 
we get\, $g+c_n u_1^2\in\iim(U)$; then\, $g\in\iim(U)$ since $u_1^2=(u_1+x)(u_1-x)$.

Substitute $g\in\iim(U)$ into (1), we get\, $f=h u_n (u_1+x)^2 u_2^2 \cdott u_{n-1}^2$\hfill (3)
\newline  for some $h$ with deg$(h)=k-1$. Then\, $f u_n=0$ gives

$h u_n^2 (u_1+x)^2 u_2^2 \cdott u_{n-1}^2=0$.

Apply $A(n,k-1)$ to $U$, CDS extra support of one, we have\, $h\in\iim(U)$.
When $k=1,\, h=0,\, f=0$. When $k\ge 2$, substitute $h\in\iim(U)$ into (3), 
$f=d  u_n^2 (u_1+x)^2 u_2^2 \cdott u_{n-1}^2$\; for some $d$.
That is, $f\in\iim_{2n-2+k}(U^{2n})$.\qqed

\lss\cbl Claim: $B(n,k)$ implies $A(n,k)$ for all $k\ge 0$.\cbk

Suppose $U$ satisfies condition (A) and $f\in \krr_k(U^{2n})$.
Choose a variable $x$ and a linear basis $\{u_1+x, u_2,u_3,\cdott, u_n\}$ of  $U$ 
such that\, $\prr_x u_i=0$  for all $1\le i\le n$. 
Observe\, $U^{2n}=$ span$((u_1+x)^2 u_2^2 \cdott u_n^2)$.
Let $g=R_{(u_1+x,\,x)} f$, then $g\in\krr_k(U^{2n})$ also.
Take  $\prr_x$ to  $g(u_1+x)^2 u_2^2 \cdott u_n^2=0$,\, we get\newline
$g u_1 u_2^2 \cdott u_n^2=0$.  So  $g u_2^2 \cdott u_n^2\in \krr_{2n-2+k}(E_x (U))$.

Since ms$\ho_1(U)\ge 2$, $x\notin U$, so dim$(E_x(U))=n$.
Apply $B(n,k)$ to $E_x (U)$, CDS exact match.
We have  $g u_2^2 \cdott u_n^2\in \iim_{2n-2+k}(E_x (U)^{2n})$.

So when $k\le 1,\, g u_2^2 \cdott u_n^2=0$ since $2n-2+k<2n$; and
when $k\ge 2$, $g u_2^2 \cdott u_n^2=h u_1^2 u_2^2 \cdott u_n^2$\; 
for some $h$, then $(g-h u_1^2)u_2^2 \cdott u_n^2=0$.

Apply $A(n-1,k)$  to  span$(u_2,\cdott, u_n)\sub U$, CDS exact match. 

When $k=0,\, g=0,\, f=0$. When $k=1,\, g\in\iim(U),\, f\in\iim(U)$. 
When\lb $k\ge 2,\, g-h u_1^2\in\iim(U)$. Since 
$u_1^2=(u_1+x)(u_1-x)$,\, $g\in\iim(U),\, f\in\iim(U)$. \qqed

\lss
Proof of Theorem 3:
Let $U$ be the linear form space spanned by the rows of  $(P\; R\; S\; T)_{n\times 4n}$,
then ms$\ho_i(U)\ge 4 i$\, for all $1\le i\le n$. By applying Theorem 7\,(A)
with $k=0$,  we have $U^{2n}\neq 0$, which gives full perrank.\qqed

\newpage\cbl
{\bf Acknowledgements}\cbk\; 

The author would like to thank Professor Noga Alon for valuable discussions and comments 
on this paper, during the 10th Krakow Conference on Graph Theory, September 2025.

The author also thanks Micha Christoph, Jeff Kahn, Noah Kravitz and Peter Pach
for helpful comments and suggestions.

\bigskip\bigskip
\centerline{\large\bf References}
\renewcommand{\refname}{} 

\end{document}